\documentclass[12pt]{article}
\usepackage{amsmath}
\usepackage{epsf}
\usepackage{amssymb}
\usepackage{graphicx}
\usepackage{ifpdf}

\newtheorem{theorem}{Theorem}[section]

\newtheorem{lemma}{Lemma}[section]

\newtheorem{property}{Property}[section]

\begin{document}
\pagestyle{empty}
\renewcommand{\thefootnote}{\fnsymbol{footnote}}

\begin{center}
{\bf \Large The generalized 3-connectivity of the folded hypercube $FQ_n$}\footnote{This
work was supported by  Hunan Education Department Foundation(No.18A382).}
\vskip 5mm

{{\bf Jing Wang$^{1}$, Fangmin Li$^{1,2}$ }\\[2mm]
$^1$ School of Computer Engineering and Applied Mathematics, Changsha University, Changsha 410022, China\\
$^2$ Hunan Province Key Laboratory of Industrial Internet Technology and Security, Changsha University, Changsha 410022, China}\\[6mm]
\end{center}
\date{}

\noindent{\bf Abstract} The generalized $k$-connectivity of a graph $G$, denoted by $\kappa_k(G)$, is a generalization of the traditional connectivity. It is well known that the generalized $k$-connectivity is an important indicator for measuring the fault tolerance and reliability of interconnection networks. The $n$-dimensional folded hypercube $FQ_n$ is obtained from the $n$-dimensional hypercube $Q_n$ by adding an edge between any pair of vertices with complementary addresses. In this paper, we show that $\kappa_3(FQ_n)=n$ for $n\ge 2$, that is, for any three vertices in $FQ_n$, there exist $n$ internally disjoint trees connecting them.

{\bf MR(2000) Subject Classification} 05C40, 05C05

\section{Introduction}
\label{secintro}

It is well known that the underlying topology of an interconnection network can be modelled by a connected graph $G=(V(G),E(G))$, where $V(G)$ is the set of processors and $E(G)$ is the set of communication links in the network. The {\it $n$-dimensional hypercube}, denoted by $Q_n$, is the graph in which each vertex is corresponding to a distinct $n$-digit binary string $x_1x_2\cdots x_n$ on the set $\{0,1\}$, and edges exist between any pair of vertices differing by a single digit in their binary representation. Since $Q_n$ posses many attractive properties, such as recursive structure, regularity, symmetry, small diameter\cite{Leighton1992,Hao2014}, the hypercube structure is a widely used interconnection model.

\begin{figure}[htbp]
\begin{minipage}[t]{0.45\linewidth}
\centering
\resizebox{0.55\textwidth}{!} {\includegraphics{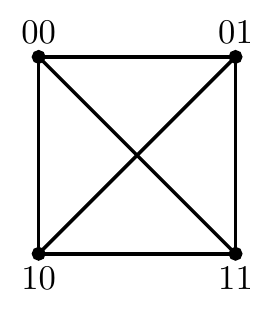}}
\caption{\small The graph $FQ_{2}$ } \label{figFQ2}
\end{minipage}
\begin{minipage}[t]{0.5\linewidth}
\centering
\resizebox{0.8\textwidth}{!} {\includegraphics{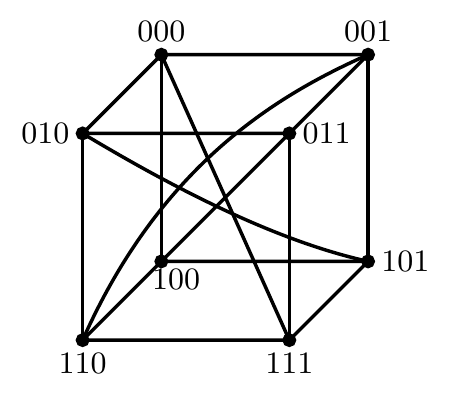}}
\caption{\small The graph $FQ_{3}$} \label{figFQ3}
\end{minipage}
\end{figure}

As a variant of the hypercube, the {\it $n$-dimensional folded hypercube} $FQ_n$, proposed first by EL-Amawy and Latifi \cite{El1991}, is a graph obtained from the hypercube $Q_n$ by adding an edge between any two vertices $x_1x_2\cdots x_n$ and $\overline{x_1} ~\overline{x_2}\cdots \overline{x_n}$, where $\overline{x_i}=1-x_i$ for all $i\in [n]$. The graphs shown in Figure \ref{figFQ2} and Figure \ref{figFQ3} are $FQ_2$ and $FQ_3$, respectively. It is shown that $FQ_n$ is $(n+1)$-regular $(n+1)$-connected. Moreover, like $Q_n$, $FQ_n$ is a Cayley graph and thus it is vertex-transitive. In \cite{El1991}, EL-Amawy and Latifi proved that $FQ_n$ has diameter $\lceil\frac{n}{2} \rceil$, about a half of the diameter of $Q_n$. Therefore, $FQ_n$ is superior to $Q_n$ in some properties.

In the design and analysis of interconnection network, the traditional connectivity of the corresponding graph is an essential parameter measuring reliability and fault tolerance of the network. For a graph $G=(V(G),E(G))$, the {\it connectivity} $\kappa(G)$ of $G$ is the minimum cardinality of a set $V'\subseteq V(G)$ such that $G-V'$ is disconnected or trivial. A famous theorem of Whitney \cite{Whitney1932} provides an equivalent definition of the connectivity. For each 2-subsets $S=\{X,Y\}$ of $V(G)$, let $\kappa_G(S)$ denote the maximum number of internally disjoint $(X,Y)$-paths in $G$. Then
$$\kappa(G)=\min\{\kappa_G(S):~ S\subseteq V(G) ~{\rm and}~ |S|=2\}.$$
As a generalization of the traditional connectivity, Chartrand et al. \cite{Chartrand1984} introduced the concept of the generalized $k$-connectivity. Let $G=(V(G),E(G))$ be a connected graph and let $k$ be an integer with $2\le k\le |V(G)|$. For a set $S$ of $k$ vertices of $G$, a tree $T$ in $G$ is called an $S$-tree, if $S\subseteq V(T)$. A collection of $\{T_1, T_2, \cdots, T_r\}$ of $S$-trees in $G$ is called internally disjoint if $E(T_i)\cap E(T_j)=\emptyset$ and $V(T_i)\cap V(T_j)=S$ for any $1\le i\ne j\le r$. Let $\kappa_G(S)$ denote the maximum number of internally disjoint $S$-trees in $G$. The {\it generalized $k$-connectivity}, denoted by $\kappa_k(G)$, is defined to be
$$\kappa_k(G)=\min\{\kappa_G(S): S\subseteq V(G) ~{\rm and}~ |S|=k\}.$$
Clearly, the generalized 2-connectivity $\kappa_2(G)$ is exactly the traditional connectivity $\kappa(G)$.

Over the past few years, the generalized $k$-connectivity has received meaningful progress. Li et al. \cite{SLi2012n} derived that it is NP-complete for a general graph $G$ to decide whether there are $l$ internally disjoint trees connecting $S$, where $l$ is a fixed integer and $S\subseteq V(G)$. Authors in \cite{HZLi2014,SLi2010} investigated the bounds of $\kappa_k(G)$ and the relationship between $\kappa(G)$ and $\kappa_k(G)$ for a general graph $G$. Many authors studied exact values of the generalized connectivity of some important classes of graphs, such as, complete graphs \cite{Chartrand2010}, complete bipartite graphs \cite{SLi2012b}, Cartesian product graphs \cite{HZLi2012,HZLi2017}, star graphs and bubble-sort graphs \cite{SLi2016}, bubble-sort star graphs \cite{Hao20191}, $(n,k)$-bubble-sort graphs \cite{Hao20192}, Cayley graphs generated by trees and cycles \cite{SLi2017}, hypercubes \cite{HZLi2012,SLin2017,Roskind1985}, dual cubes \cite{ZhaoHao2019}, and so on.

In this paper, we focus on determining the generalized 3-connectivity of the folded hypercube $FQ_n$, and get the main result as follows.

\begin{theorem}\label{thk3FQn}
$\kappa_3(FQ_n)=n$ for any $n\ge 2$.
\end{theorem}

\section{Preliminaries}\label{secpreli}
Throughout this paper, we consider a simple, connected graph $G=(V(G),E(G))$ with $V(G)$ be its vertex set and $E(G)$ be its edge set. For a vertex $X\in V(G)$, the {\it degree} of $X$ in $G$, denoted by ${\rm deg}_G(X)$, is the number of edges of $G$ incident with $X$. Let $V'\subseteq V(G)$, denote by $G-V'$ the graph obtained from $G$ by deleting all the vertices in $V'$ together with their incident edges. Let $P=X_0X_1\cdots X_m$ be a path connecting two vertices $X_0$ and $X_m$ in $G$, then $P$ is called an $(X_0,X_m)$-{\it path}.

Let $n\ge 2$ and $d$ be an integer that $1\le d\le n$. Denote $Q_n^d[0]$ (resp. $Q_n^d[1]$) be the graph with vertex set $V(Q_n^d[0])=\{x_1\cdots x_{d-1}0x_{d+1}\cdots x_n: x_i\in\{0,1\}, i\in [n]\backslash \{d\}\}$ (resp. $V(Q_n^d[1])=\{x_1\cdots x_{d-1}1x_{d+1}\cdots x_n: x_i\in\{0,1\}, i\in [n]\backslash \{d\}\}$), and two vertices in $Q_n^d[0]$ (resp. $Q_n^d[1]$) are adjacent if and only if they differ in exactly one coordinate. It is easily seen that Property \ref{proQnd} holds.

\begin{property}\label{proQnd}
(1) Both $Q_n^d[0]$ and $Q_n^d[1]$ are induced subgraphs of the hypercube $Q_n$;\\
(2) Both $Q_n^d[0]$ and $Q_n^d[1]$ are isomorphic to $Q_{n-1}$.
\end{property}

Let $X=x_1x_2\cdots x_n$ and $Y=y_1y_2\cdots y_n$ be two vertices in the hypercube $Q_n$. If $X$ and $Y$ exactly differ by the $d$th digit ($1\le d\le n$), that is $|x_d-y_d|=1$ and $x_i=y_i$ for all $i\in [n]\backslash \{d\}$, then $XY$ is named a {\it $d$-dimensional edge}. Based on these notations, we can divide $Q_n$ along the $d$th dimension into $Q_n^d[0]$ and $Q_n^d[1]$ by deleting all the $d$-dimensional edges from $Q_n$.

The folded hypercube $FQ_n$ ($n\ge 2$) can be obtained from $Q_n^d[0]$ and $Q_n^d[1]$ by adding $2\times 2^{n-1}$ edges as follows: a vertex $X=x_1\cdots x_{d-1}0x_{d+1}\cdots x_n\in V(Q_n^d[0])$ is joined to a vertex $Y=y_1\cdots y_{d-1}1y_{d+1}\cdots y_n\in V(Q_n^d[1])$ if for every $i\in [n]\backslash \{d\}$, either\\
$(i)$~ $x_i=y_i$; in this case, $XY$ is called a {\it hypercube edge}, or\\
$(ii)$~ $x_i=\overline{y_i}$; in this case, $XY$ is called a {\it complement edge}.\\
For any integer $d\in [n]$, the construction of $FQ_n$ can be symbolically written as $FQ_n=Q_n^d[0]\bigotimes Q_n^d[1]$. It follows from the definition that $X$ differs with $Y$ by $d$th digit when $XY$ is a hypercube edge, in this case, denote $Y=N_d(X)$. Furthermore, $X$ differs with $Y$ by each digit when $XY$ is a complement edge, in this case, denote $Y=\overline{X}$.

%
%

Let $X$ be an arbitrary vertex in $V(FQ_n)\cap V(Q_n^d[0])$, it has two neighbours in $Q_n^d[1]$, one is $N_d(X)$, which is called a {\it hypercube neighbour} of $X$, the other is $\overline{X}$, which is called a {\it complement neighbour} of $X$ . Both $N_d(X)$ and $\overline{X}$ are named {\it out neighbours} of $X$ in $Q_n^d[1]$. Similarly, for any vertex $Y$ in $V(FQ_n)\cap V(Q_n^d[1])$, it has a hypercube neighbour $N_d(Y)$ and a complement neighbour $\overline{Y}$ in $Q_n^d[0]$.

The authors in \cite{Bondy,HZLi2012} studied the connectivity and the generalized 3-connectivity of the hypercube $Q_n$, respectively.

\begin{lemma}\label{leQn2}(\cite{Bondy})
For any integer $n\ge 2$, $\kappa(Q_n)=n$.
\end{lemma}

\begin{lemma}\label{leQn3}(\cite{HZLi2012})
For any integer $n\ge 2$, $\kappa_3(Q_n)=n-1$.
\end{lemma}

Li et al. \cite{SLi2010} gave an upper bound of $\kappa_3(G)$ for a general graph $G$.

\begin{lemma}\label{leupperK}(\cite{SLi2010})
Let $G$ be a connected graph with minimum degree $\delta$. If there are two adjacent vertices of degree $\delta$, then $\kappa_3(G)\le \delta-1$.
\end{lemma}

Chartrand et al. \cite{Chartrand2010} established the generalized $k$-connectivity of the complete graph $K_n$.

\begin{lemma}(\cite{Chartrand2010})\label{leKn}
For every two integers $n$ and $k$ with $2\le k\le n$, $\kappa_k(K_n)=n-\lceil\frac{k}{2}\rceil$.
\end{lemma}

Figure \ref{figFQ2} shows that $FQ_2$ is isomorphic to the complete graph $K_4$. By Lemma \ref{leKn}, we have

\begin{lemma}\label{leFQ2}
$\kappa_3(FQ_2)=2$.
\end{lemma}

\begin{figure}[htbp]
\begin{minipage}[t]{0.5\linewidth}
\centering
\resizebox{0.8\textwidth}{!} {\includegraphics{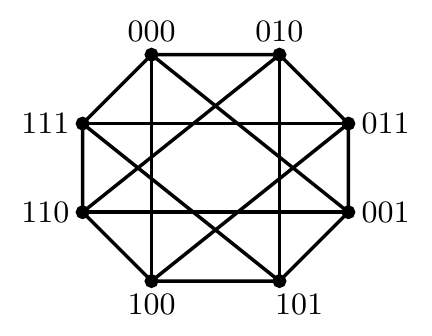}}
\caption{\small A symmetric drawing of $FQ_{3}$ } \label{figAnFQ3}
\end{minipage}
\begin{minipage}[t]{0.5\linewidth}
\centering
\resizebox{0.8\textwidth}{!} {\includegraphics{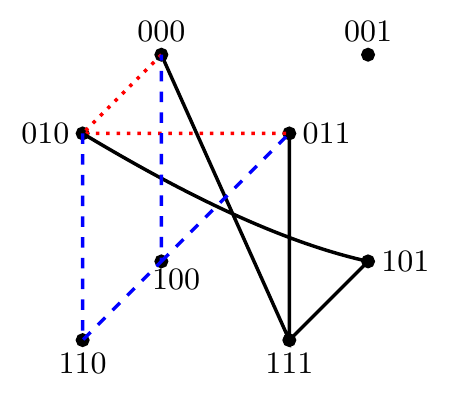}}
\caption{\small The desired $\{000,010,011\}$-trees} \label{figS1}
\end{minipage}
\end{figure}

\begin{figure}[htbp]
\begin{minipage}[t]{0.5\linewidth}
\centering
\resizebox{0.8\textwidth}{!} {\includegraphics{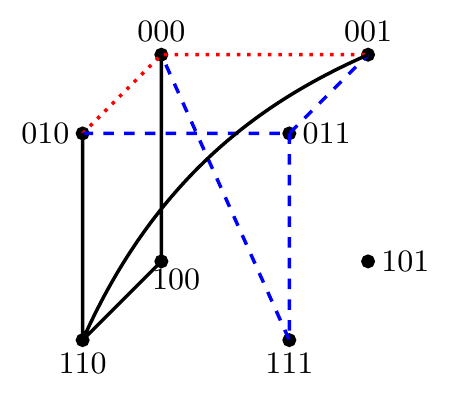}}
\caption{\small The desired $\{000,010,001\}$-trees} \label{figS2}
\end{minipage}
\begin{minipage}[t]{0.5\linewidth}
\centering
\resizebox{0.8\textwidth}{!} {\includegraphics{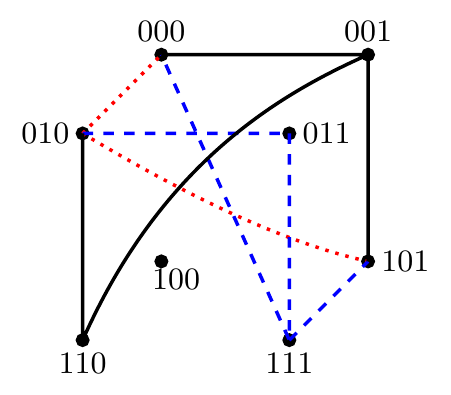}}
\caption{\small The desired $\{000,010,101\}$-trees} \label{figS3}
\end{minipage}
\end{figure}

We end this section with the following lemma that determines the lower bound of $\kappa_3(FQ_3)$.

\begin{lemma}\label{leFQ3}
Let $S$ be a subset of $V(FQ_3)$ with $|S|=3$. Then there exist three internally disjoint $S$-trees in $FQ_3$.
\end{lemma}

\noindent{\bf Proof.} Figure \ref{figAnFQ3} shows a symmetric drawing of $FQ_3$. By the symmetry of $FQ_3$, it suffices to give three internally disjoint $S$-trees for each of the following 3-subsets $S\subseteq V(FQ_3)$:
$$\{000,010,011\},\; \{000,010,001\}, \; \{000,010,101\}, \; \{000,011,101\}, \; \{000,011,100\}.$$

\vskip 1mm

\begin{figure}[htbp]
\begin{minipage}[t]{0.5\linewidth}
\centering
\resizebox{0.8\textwidth}{!} {\includegraphics{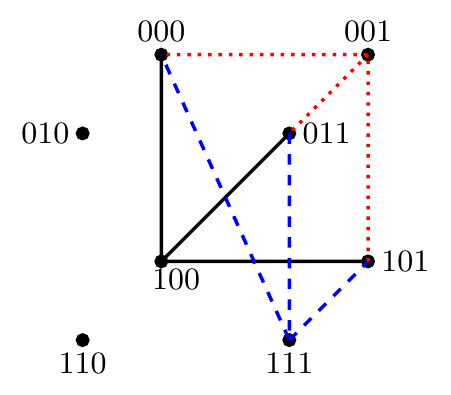}}
\caption{\small The desired $\{000,011,101\}$-trees} \label{figS4}
\end{minipage}
\begin{minipage}[t]{0.5\linewidth}
\centering
\resizebox{0.8\textwidth}{!} {\includegraphics{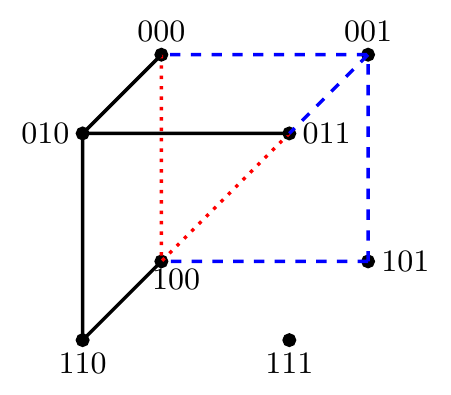}}
\caption{\small The desired $\{000,011,100\}$-trees} \label{figS5}
\end{minipage}
\end{figure}

If $S=\{000,010,011\}$, then three internally disjoint $S$-trees are depicted in Figure \ref{figS1}, where three trees are represented by black solid lines, blue dashed lines and red dotted lines, respectively. If $S=\{000,010,001\}$, $S=\{000,010,101\}$, $S=\{000,011,101\}$ or $S=\{000,011,100\}$, the desired $S$-trees are shown in Figures \ref{figS2}, \ref{figS3}, \ref{figS4} or \ref{figS5}, respectively.   \hfill$\Box$

\section{Proof of Theorem \ref{thk3FQn}}

Denote $Q_n^{12}[10]$ be the graph whose vertex set is $\{10x_3\cdots x_n: x_i\in\{0,1\}, 3\le i\le n\}$, where two vertices are adjacent if and only if they differ in precisely one coordinate. Similarly, let $Q_n^{12}[11]$ be the graph whose vertex set is $\{11x_3\cdots x_n: x_i\in\{0,1\}, 3\le i\le n\}$, where two vertices are adjacent if and only if they differ in precisely one coordinate. Obviously, both $Q_n^{12}[10]$ and $Q_n^{12}[11]$ are subgraphs of $Q_n^1[1]$ and are isomorphic to $Q_{n-2}$.

We define some special binary string vertex which will be used in latter discussions. For every integer $1\le i\le n$, let $E_{n,i}$ (or $E_i$, for simplicity) denote the binary string vertex of $FQ_n$ with the $i$th digit being 1 and all the remaining digits being 0. That is,
$$E_1=100\cdots 00, \; E_2=010\cdots 00, \; E_n=000\cdots 01.$$
Furthermore, let $W_i=\overline{E_i}$ for all $1\le i\le n$, which is the complement neighbour of $E_i$ in $FQ_n$. That is,
$$W_1=011\cdots 11, \; W_2=101\cdots 11, \; W_n=111\cdots 10.$$

The following lemma is dedicated to give the lower bound of $\kappa_3(FQ_n)$ for $n\ge 4$.

\begin{lemma}\label{lelowerFQn}
For $n\ge 4$, let $S=\{X,Y,Z\}$ be any 3-subsets of $V(FQ_n)$. There exist $n$-internally disjoint $S$-trees in $FQ_n$.
\end{lemma}

\noindent{\bf Proof}. Note that $FQ_n=Q_n^d[0]\bigotimes Q_n^d[1]$, $d\in [n]$. By the symmetry of the folded hypercube, we consider that $FQ_n=Q_n^1[0]\bigotimes Q_n^1[1]$. W.L.O.G, assume that $|S\cap V(Q_n^1[0])|\ge |S\cap V(Q_n^1[1])|$. The following two cases are considered.

\vskip 2mm
{\bf Case 1}. $|S\cap V(Q_n^1[0])|=3$.

That is, $X=0x_2\cdots x_n\triangleq 0A$, $Y=0y_2\cdots y_n\triangleq 0B$ and $Z=0z_2\cdots z_n\triangleq 0C$, where $A=x_2\cdots x_n$,$B=y_2\cdots y_n$ and $C=z_2\cdots z_n$ are $(n-1)$-digit binary strings. Combined with Property \ref{proQnd} and Lemma \ref{leQn3}, there exist ($n-2$) internally disjoint $S$-trees $T_1, T_2, \cdots, T_{n-2}$ in $Q_n^1[0]$.

Now we are going to find two other $S$-trees in $FQ_n$. Recall that $X=0A$ has two neighbours $N_1(X)=1A$ and $\overline{X}=1\overline{A}$ in $Q_n^1[1]$, $Y=0B$ (resp. $Z=0C$) has two neighbours $N_1(Y)=1B$ and $\overline{Y}=1\overline{B}$ (resp. $N_1(Z)=1C$ and $\overline{Z}=1\overline{C}$) in $Q_n^1[1]$. Moreover, $Q_n^1[1]$ can be obtained from two disjoint copies of $Q_n^{12}[10]$ and $Q_n^{12}[11]$ by adding a perfect matching. W.L.O.G., assume that the first binary digit of $A, B$ and $C$ is 0. That is,
$$\{N_1(X),N_1(Y),N_1(Z)\}\subseteq V(Q_n^{12}[10]), \;\; {\rm and}\;\; \{\overline{X},\overline{Y},\overline{Z}\}\subseteq V(Q_n^{12}[11]).$$

There exists an $\{N_1(X),N_1(Y),N_1(Z)\}$-tree $\hat{T}_{n-1}$ in $Q_n^{12}[10]$ since $Q_n^{12}[10]$ is isomorphic to $Q_{n-2}$. Analogously, there exists an $\{\overline{X},\overline{Y},\overline{Z}\}$-tree $\hat{T}_{n}$ in $Q_n^{12}[11]$. Let $T_{n-1}=\hat{T}_{n-1}+\{XN_1(X), YN_1(Y), ZN_1(Z)\}$ and $T_{n}=\hat{T}_{n}+\{X\overline{X}, Y\overline{Y}, Z\overline{Z}\}$. It is not difficult to testify that $T_1, \cdots, T_{n-1}, T_{n}$ are internally disjoint $S$-trees.

\vskip 2mm
{\bf Case 2}. $|S\cap V(Q_n^1[0])|=2$.

W.L.O.G., assume that $\{X,Y\}\subseteq V(Q_n^1[0])$ and $Z\in V(Q_n^1[1])$. Furthermore, let $Z=11\cdots 11$, $X=0x_2\cdots x_n$ and $Y=0y_2\cdots y_n$.

Remind that $Z$ has two neighbours, $N_1(Z)=01\cdots 11$ and $\overline{Z}=00\cdots 00$, in $Q_n^1[0]$.

\vskip 2mm
{\bf Subcase 2.1}. $|\{N_1(Z), \overline{Z}\}\cap \{X,Y\}|=2$.

That is $\{N_1(Z), \overline{Z}\}=\{X,Y\}$. W.L.O.G., assume that $X=N_1(Z)$ and $Y=\overline{Z}$. One can see that $X$ is not adjacent to $Y$ since $n\ge 4$. By Lemma \ref{leQn2}, there exist ($n-1$) internally disjoint ($X,Y$)-paths $P_2, \cdots, P_{n}$ in $Q_n^1[0]$. Let $N_2(X), N_3(X), \cdots, N_{n}(X)$ be the neighbours of $X$ such that $N_i(X)$ is in $P_i$. Note that $W_i$, the hypercube neighbour of $N_i(X)$ in $Q_n^1[1]$, is adjacent to $Z$.

Let $T_1=\{XZ, YZ\}$ and let $T_i=P_i+\{N_i(X)W_i, W_iZ\}$ for $2\le i\le n$. Clearly, $T_1, T_{2}, \cdots, T_{n}$ are internally disjoint $S$-trees.

\vskip 2mm
{\bf Subcase 2.2}. $|\{N_1(Z), \overline{Z}\}\cap \{X,Y\}|=1$.

\vskip 2mm
{\bf Subcase 2.2.1}. $\{N_1(Z), \overline{Z}\}\cap \{X,Y\}=X=N_1(Z)$.

By similarly arguments in Subcase 2.1, there exist ($n-1$) internally disjoint ($X,Y$)-paths $P_2, \cdots, P_{n}$ in $Q_n^1[0]$. Let $N_2(X), N_3(X), \cdots, N_{n}(X)$ be the neighbours of $X$ in $Q^1_n[0]$ and $N_i(X)$ is in $P_i$. Note that $N_1(N_i(X))=W_i$ is adjacent to $Z$.

It is possible that $X$ is adjacent to $Y$, we may assume $Y=N_2(X)$ under this circumstance. Let $T_2=P_2+\{XZ\}$, and let $T_i=P_i+\{N_i(X)W_i, W_iZ\}$ for $3\le i\le n$.

Set $W'=\bigcup_{i=3}^{n}\{W_i\}$. It is observed that $\overline{X}\notin W'$, and at least one of $N_1(Y)$ and $\overline{Y}$ is not contained in $W'$ since $n\ge 4$. W.L.O.G., assume that $\overline{Y}\notin W'$. According to Lemma \ref{leQn2}, $Q^1_n[1]-W'$ is connected since $|W'|=n-2$. Consequently, there is an $\{\overline{X}, \overline{Y}, Z\}$-tree $\hat{T}_{1}$ in $Q^1_n[1]-W'$. Let $T_1=\hat{T}_{1}+\{X\overline{X}, Y\overline{Y}\}$. Therefore, $T_1, T_{2}, \cdots, T_{n}$ are internally disjoint $S$-trees.

\vskip 2mm
{\bf Subcase 2.2.2}. $\{N_1(Z), \overline{Z}\}\cap \{X,Y\}=X=\overline{Z}$.

That is $X=00\cdots 00$, and $E_2,\cdots, E_n$ are ($n-1$) neighbours of $X$ in $Q^1_n[0]$. It is possible that $X$ is adjacent to $Y$, we may assume $Y=E_2$ under this circumstance. Let $W'=\bigcup_{i=3}^{n}\{\overline{E_i}\}=\bigcup_{i=3}^{n}\{W_i\}$. By letting $T_i=P_i+\{E_iW_i, W_iZ\}$ for $3\le i\le n$ and constructing $T_1$ and $T_2$ as in Subcase 2.2.1, $n$ internally disjoint $S$-trees can be obtained.

\vskip 2mm
{\bf Case 2.3}. $|\{N_1(Z), \overline{Z}\}\cap \{X,Y\}|=0$.

That is
\begin{eqnarray}\label{eq1}
\{N_1(Z), \overline{Z}\}\cap \{X,Y\}=\emptyset.
\end{eqnarray}

Note that $N_2(Z)=W_2, N_3(Z)=W_3, \cdots, N_n(Z)=W_n$ are ($n-1$) neighbours of $Z$ in $Q^1_n[1]$. Set
\begin{eqnarray}\label{eqW}
W=\bigcup_{i=2}^{n}\{W_i\}.
\end{eqnarray}

Consider two out neighbours of $X=0x_2\cdots x_n$ in $Q^1_n[1]$, $N_1(X)=1x_2\cdots x_n$ and $\overline{X}=1\overline{x_2}\cdots \overline{x_n}$. It is obviously that $|\{N_1(X), \overline{X}\}\cap W|\le 1$. Similarly, $|\{N_1(Y), \overline{Y}\}\cap W|\le 1$. For convenience, we assume that $|\{N_1(X), \overline{X}\}\cap W|\ge |\{N_1(Y), \overline{Y}\}\cap W|$.

\vskip 2mm
{\bf Subcase 2.3.1}. $|\{N_1(X), \overline{X}\}\cap W|=0$.

That is to say, $\{N_1(X), \overline{X}\}\cap W=\{N_1(Y), \overline{Y}\}\cap W=\emptyset$. Recall that $N_1(Z)=01\cdots 11$, and $N_1(W_2), \cdots, N_1(W_n)$ are ($n-1$) neighbours of $N_1(Z)$ in $Q_n^1[0]$.

By Lemma \ref{leQn3}, there are ($n-2$) internally disjoint $\{X,Y,N_1(Z)\}$-trees $\hat{T}_3, \cdots, \hat{T}_n$ in $Q_n^1[0]$. Since ${\rm deg}_{Q_n^1[0]}(N_1(Z))=n-1$, we have that $1\le {\rm deg}_{\hat{T}_i}(N_1(Z))\le 2$ for any $i\in \{3,\cdots, n\}$. 

\vskip 2mm
{\bf Subcase 2.3.1.1}. ${\rm deg}_{\hat{T}_i}(N_1(Z))=1$ for any $i\in \{3,\cdots, n\}$.

There exists a neighbour of $N_1(Z)$ in $Q_n^1[0]$, say $N_1(W_2)$, not in any tree $\hat{T}_3, \cdots, \hat{T}_n$. W.L.O.G., assume that $N_1(W_i)$ be a neighbour of $N_1(Z)$ in $\hat{T}_i$, $i\in \{3,\cdots, n\}$. Let $T_3=\hat{T}_3+\{N_1(Z)Z\}$, and let $T_i=(\hat{T}_i-N_1(Z))+\{N_1(W_i)W_i, W_iZ\}$ for $4\le i\le n$.

We may divide $Q_n^1[1]$ along the 2nd dimension into $Q_n^{12}[10]$ and $Q_n^{12}[11]$. It is seen that $W_2\in V(Q_n^{12}[10])$ and $W_3\in V(Q_n^{12}[11])$. Now consider $N_1(X)$ and $\overline{X}$, two out neighbours of $X$ in $Q_n^1[1]$. They belong to different copies of $Q_n^{12}[10]$ and $Q_n^{12}[11]$. W.L.O.G., assume that  $N_1(X)\in V(Q_n^{12}[10])$ and $\overline{X}\in Q_n^{12}[11]$. Analogously, let $N_1(Y)\in V(Q_n^{12}[11])$ and $\overline{Y}\in Q_n^{12}[10]$. It is possible that $N_1(X)=\overline{Y}$.

Let $$H=\bigcup_{i=4}^{n}\{W_i\}.$$
Clearly, $H\subseteq V(Q_n^{12}[11])$ and $|H|=n-3$.

There exists an $\{N_1(X), \overline{Y}, W_2\}$-tree $\hat{T}_{1}$ in $Q_n^{12}[10]$ avoiding vertices of $H$ since $H\cap V(Q_n^{12}[10])=\emptyset$. Let $T_1=\hat{T}_{1}+\{XN_1(X), Y\overline{Y}, W_2Z\}$.

According to the fact that $Q_n^{12}[11]$ is isomorphic to $Q_{n-2}$ and by Lemma \ref{leQn2}, we know that $Q_n^{12}[11]-H$ is still connected. Hence, there exists an $\{\overline{X}, N_1(Y), Z\}$-tree $\hat{T}_{2}$ in $Q_n^{12}[11]-H$. Let $T_2=\hat{T_2}+\{X\overline{X}, YN_1(Y)\}$.

Clearly, $T_1, T_2, \cdots, T_n$ are $n$ internally disjoint $S$-trees.

\vskip 2mm
{\bf Subcase 2.3.1.2}. There exists an integer $i\in\{3,\cdots, n\}$ such that ${\rm deg}_{\hat{T}_i}(N_1(Z))=2$.

W.L.O.G., assume that ${\rm deg}_{\hat{T}_3}(N_1(Z))=2$ and $N_1(W_2)$ and $N_1(W_3)$ be two neighbours of $N_1(Z)$ in $\hat{T}_3$. Then ${\rm deg}_{\hat{T}_i}(N_1(Z))=1$ for any $4\le i\le n$.

Define $T_i$ as that in Subcase 2.3.1.1, $n$ internally disjoint $S$-trees $T_1, T_2, \cdots, T_n$ can be obtained.

\vskip 2mm
{\bf Subcase 2.3.2}. $|\{N_1(X), \overline{X}\}\cap W|=1$ and $|\{N_1(Y), \overline{Y}\}\cap W|=1$.

\vskip 2mm
{\bf Subcase 2.3.2.1}. $\{N_1(X), \overline{X}\}\cap W=N_1(X)$ and $\{N_1(Y), \overline{Y}\}\cap W=N_1(Y)$.

W.L.O.G., assume that $N_1(X)=W_2=1011\cdots 1$ and $N_1(Y)=W_3=1101\cdots 1$. That means $X=0011\cdots 1$ and $Y=0101\cdots 1$. Since $n\ge 4$, we have that the $n$th digit of $X, Y$ and $Z$ is 1, respectively. Then we may divide $FQ_n$ along the $n$th dimension into $Q_n^n[0]$ and $Q_n^n[1]$ such that $|S\cap V(Q_n^n[1])|=3$. Similarly to Case 1, we can find $n$ internally disjoint $S$-trees.

\vskip 2mm
{\bf Subcase 2.3.2.2}. $\{N_1(X), \overline{X}\}\cap W=N_1(X)$ and $\{N_1(Y), \overline{Y}\}\cap W=\overline{Y}$.

W.L.O.G., assume that $N_1(X)=W_2=1011\cdots 1$ and $\overline{Y}=W_3=1101\cdots 1$. That is, $X=0011\cdots 1$ and $Y=0010\cdots 0$. Then we may divide $FQ_n$ along the 3rd dimension into $Q_n^3[0]$ and $Q_n^3[1]$ such that $|S\cap V(Q_n^3[1])|=3$. Similarly to Case 1, we can find $n$ internally disjoint $S$-trees.

\vskip 2mm
{\bf Subcase 2.3.2.3}. $\{N_1(X), \overline{X}\}\cap W=\overline{X}$ and $\{N_1(Y), \overline{Y}\}\cap W=\overline{Y}$.

W.L.O.G., assume that $\overline{X}=W_2=1011\cdots 1$ and $\overline{Y}=W_3=1101\cdots 1$. That is, $X=0100\cdots 0=E_2$ and $Y=0010\cdots 0=E_3$.

Firstly, we can construct ($n-1$) internally disjoint ($X,Y$)-paths $P_2, P_3, \cdots, P_n$ in $Q_n^1[0]$ as follows:\\
$P_2$: $(X=01000\cdots 00, 00000\cdots 00=\overline{Z}, 00100\cdots 00=Y)$;\\
$P_3$: $(X=01000\cdots 00, 01100\cdots 00, 00100\cdots 00=Y)$;\\
$P_4$: $(X=01000\cdots 00, 01010\cdots 00, 00010\cdots 00, 00110\cdots 00, 00100\cdots 00=Y)$;\\
$\cdots \cdots$\\
$P_n$: $(X=01000\cdots 00, 01000\cdots 01, 00000\cdots 01, 00100\cdots 01, 00100\cdots 00=Y)$.

For $4\le i\le n$, we have that $P_i=(X,N_i(X),N_2(N_i(X)),N_3(N_2(N_i(X))),Y)$. Note that $N_2(N_i(X))=E_i$ and $\overline{E_i}=W_i$. Based on the paths $P_2, P_3, \cdots, P_n$, we can define $T_2=P_2+\{\overline{Z}Z\}$, $T_3=P_3+\{YW_3, W_3Z\}$, and $T_i=P_i+\{E_iW_i, W_iZ\}$ for $4\le i\le n$.

Let $W'=\bigcup_{i=3}^{n}\{W_i\}$. It is obviously that $|W'|=n-2$, $N_1(X)\notin W'$ and $N_1(Y)\notin W'$. By Lemma \ref{leQn2}, $Q_n^1[1]-W'$ is connected since $\kappa(Q_n^1[1])=\kappa(Q_{n-1})=n-1$. Thus, there exists an $\{N_1(X), N_1(Y), Z\}$-tree $\hat{T}_1$ in $Q_n^1[1]-W'$. Let $T_1=\hat{T}_1+\{XN_1(X), YN_1(Y)\}$. Clearly, $T_1, T_2, \cdots, T_n$ are $n$ internally disjoint $S$-trees.

\vskip 2mm
{\bf Subcase 2.3.2.4}. $\{N_1(X), \overline{X}\}\cap W=\overline{X}$ and $\{N_1(Y), \overline{Y}\}\cap W=N_1(Y)$.

Under this restriction, the argument is similar to that in Subcase 2.3.2.2 and is omitted.

\vskip 2mm
{\bf Subcase 2.3.3}. $|\{N_1(X), \overline{X}\}\cap W|=1$ and $|\{N_1(Y), \overline{Y}\}\cap W|=0$.

That is to say,
\begin{eqnarray}\label{eq2}
\{N_1(Y), \overline{Y}\}\cap W=\emptyset.
\end{eqnarray}

The following two subcases are distinguished.

\vskip 2mm
{\bf Subcase 2.3.3.1}. $\{N_1(X), \overline{X}\}\cap W=N_1(X)$.

W.L.O.G., assume that $N_1(X)=W_2=101\cdots 1$. That is $X=001\cdots 1$.

Recall that $Y=0y_2y_3\cdots y_n$. If there exists an integer $i\in \{3,\cdots, n\}$ such that $y_i=1$, then we may divide $FQ_n$ along the $i$th dimension such that $|\{X,Y,Z\}\cap V(Q_n^i[1])|=3$. Similarly to Case 1, we can find $n$ internally disjoint $S$-trees. Therefore, it suffice to consider the case that $y_3=\cdots=y_n=0$, that is $Y=0y_20\cdots 0$. Furthermore, we can conclude that $y_2=1$, otherwise $Y=\overline{Z}$, a contradiction with Eq.(\ref{eq1}). Under this situation, $Y=010\cdots 0$ and $\overline{Y}=W_2\in W$, a contradiction with Eq.(\ref{eq2}).

\vskip 2mm
{\bf Subcase 2.3.3.2}. $\{N_1(X), \overline{X}\}\cap W=\overline{X}$.

W.L.O.G., assume that $\overline{X}=W_2=101\cdots 1$. That is to say, $X=010\cdots 0$. Now consider $Y=0y_2y_3\cdots y_n$.

If $y_2=1$, then we may divide $FQ_n$ along the 2nd dimension such that $|\{X,Y,Z\}\cap V(Q_n^2[1])|=3$. Similarly to Case 1, we can find $n$ internally disjoint $S$-trees.

We now consider the case that $y_2=0$, that is $Y=00y_3\cdots y_n$.

\vskip 2mm
{\bf Claim 1}. At least two among $y_3,\cdots, y_n$ are digit 1.

{\bf Proof of Claim 1}. Suppose to contrary that there is at most one among $y_3,\cdots, y_n$ is digit 1. Then $Y=\overline{Z}$ if $y_3=\cdots=y_n=0$, a contradiction with Eq.(\ref{eq1}); or $\overline{Y}=W_j$ if $y_j=1$ and $y_i=0$ for $3\le i\ne j \le n$, a contradiction with Eq.(\ref{eq2}). \hfill$\Box$

\vskip 2mm
{\bf Claim 2}. At least one among $y_3,\cdots, y_n$ is digit 0.

{\bf Proof of Claim 2}. Otherwise, $Y=001\cdots 1$ and $N_1(Y)$, the hypercube neighbour of $Y$ in $Q_n^1[1]$, is $W_2$, a contradiction with Eq.(\ref{eq2}).  \hfill$\Box$

\vskip 2mm

We may divide $FQ_n$ along the 2nd dimension such that $\{X,Z\}\subseteq V(Q_n^2[1])$ and $Y\in V(Q_n^2[0])$. The outside neighbours of $Y$ in $Q_n^2[1]$ are $N_2(Y)=01y_3\cdots y_n$ and $\overline{Y}=11\overline{y_3}\cdots \overline{y_n}$. Moreover, ($n-1$) neighbours of $Y$ in $Q_n^2[0]$ are
$$N_1(Y)=10y_3\cdots y_n,\; N_3(Y)=00\overline{y_3}\cdots y_n, \; N_n(Y)=00y_3\cdots y_{n-1}\overline{y_n}.$$

It is noted that $X$ and $Z$ have two common out neighbours in $Q_n^2[0]$, $\overline{X}=101\cdots 1$ and $\overline{Z}=000\cdots 0$. According to Claim 1 and Claim 2, it is inferred that $\{N_2(Y), \overline{Y}\}\cap \{X,Z\}=\emptyset$, and $\{\overline{Z},\overline{X}\}\cap \{N_1(Y), N_3(Y), \cdots, N_n(Y)\}=\emptyset.$ Using similar analysis in Subcase 2.3.1, $n$ internally disjoint $S$-trees can be obtained.

In all, there exist $n$ internally disjoint $S$-trees in $FQ_n$ for any 3-subsets $S\subseteq V(FQ_n)$. \hfill$\Box$

\vskip 2mm
Now, we are in a position to prove Theorem \ref{thk3FQn}.

\vskip 2mm
\noindent{\bf Proof of Theorem \ref{thk3FQn}}. Lemma \ref{leFQ2} implies that the theorem is true for $n=2$. Therefore, we only need to prove the result for $n\ge 3$.

It is known from Lemma \ref{leupperK} that $\kappa_3(FQ_n)\le \delta -1=n$ since $FQ_n$ is ($n+1$)-regular. On the other side, Lemma \ref{leFQ3} together with Lemma \ref{lelowerFQn} confirm that $\kappa_3(FQ_n)\ge n$ for $n\ge 3$. Consequently, the proof is done. \hfill$\Box$

\section{Conclusion}

The generalized $k$-connectivity is a natural generalization of the traditional connectivity and can serve for measuring the capability of a network $G$ to connect any $k$ vertices in $G$. In this paper, we focus on the generalized 3-connectivity of the folded hypercube $FQ_n$, and get the main result that $\kappa_3(FQ_n)=n$. So far, the generalized $k$-connectivity of the network are almost about $k=3$. It would be an interesting topic to study $\kappa_k(FQ_n)$ for $k\ge 4$.

%



\end{document}